\documentclass[oneside, a4paper,reqno]{amsart}
\usepackage{pdfsync}
\usepackage{stmaryrd}
\usepackage{mathrsfs}
\usepackage{multicol}
\usepackage{amsmath, amsthm, amscd, amssymb, latexsym, eucal}
\usepackage[all]{xy}
\def\serieslogo@{} \def\@setcopyright{} \makeatother

\usepackage{multienum}

\usepackage{hyperref}
\usepackage{color}
\usepackage{cite}

\hypersetup{colorlinks=true,
     breaklinks=true,
     linkcolor=blue,    
     bookmarks=true,
      pageanchor=true}

\makeatletter
\renewcommand*\env@matrix[1][c]{\hskip -\arraycolsep
  \let\@ifnextchar\new@ifnextchar
  \array{*\c@MaxMatrixCols #1}}
\makeatother

\usepackage{color}

\usepackage[colorinlistoftodos]{todonotes}

\numberwithin{equation}{section}
\newtheorem{thm}{Theorem}[section]
\newtheorem{cor}[thm]{Corollary}
\newtheorem{lem}[thm]{Lemma}

\theoremstyle{definition}

\newtheorem{rem}[thm]{Remark}

\newcommand{\lxr}{\longrightarrow}

\newcommand{\A}{\mathscr A}
\newcommand{\B}{\mathscr B}
\newcommand{\C}{\mathscr C}

\newcommand{\F}{\mathcal F}

\newcommand{\mt}{\mathsf{T}}
\newcommand{\mh}{\mathsf{H}}
\newcommand{\mU}{\mathsf{U}}

\newcommand{\mr}{\mathsf{r}}
\newcommand{\mq}{\mathsf{q}}
\newcommand{\mi}{\mathsf{i}}

\newcommand{\ml}{\mathsf{l}}
\newcommand{\me}{\mathsf{e}}

\newcommand{\map}{\mathsf{p}}

 \DeclareMathOperator{\inc}{\mathsf{inc}}
  
\DeclareMathOperator*{\Ker}{\mathsf{Ker}}
 \DeclareMathOperator*{\Image}{\mathsf{Im}}
\DeclareMathOperator*{\Coker}{\mathsf{Coker}}

 \DeclareMathOperator{\pd}{\mathsf{pd}}

  \DeclareMathOperator*{\gld}{\mathsf{gl.dim}}

 \DeclareMathOperator*{\smod}{\mathsf{mod}-\!}
 \DeclareMathOperator*{\lsmod}{\!-\mathsf{mod}}
\DeclareMathOperator*{\lMod}{\!-\mathsf{Mod}}

\DeclareMathOperator*{\Add}{\mathsf{Add}}
 
\DeclareMathOperator*{\add}{\mathsf{add}}

\DeclareMathOperator{\Hom}{\mathsf{Hom}}
 
 \DeclareMathOperator{\rad}{\mathsf{rad}}

\DeclareMathOperator*{\Tor}{\mathsf{Tor}}

  \DeclareMathOperator*{\op}{\mathsf{op}}

\newcommand{\iden}{\operatorname{Id}\nolimits}

\newsavebox{\proofbox}
\savebox{\proofbox}{\begin{picture}(7,7)%
  \put(0,0){\framebox(7,7){}}\end{picture}}

\begin{document}

\title[]{Recollements of abelian categories and ideals in heredity
chains - a recursive approach to quasi-hereditary algebras}

\author[]{Nan Gao, Steffen Koenig and Chrysostomos Psaroudakis}
\address{Nan Gao \\ Department of Mathematics, Shanghai University, Shanghai 200444, PR China}
\email{nangao@shu.edu.cn}

\address{Steffen Koenig \\ Institute of Algebra and Number Theory, University of Stuttgart, Pfaffenwaldring 57, 70569 Stuttgart, Germany}
\email{skoenig@mathematik.uni-stuttgart.de}

\address{Chrysostomos Psaroudakis\\ Institute of Algebra and Number Theory, University of Stuttgart, Pfaffenwaldring 57, 70569 Stuttgart, Germany}
\email{chrysostomos.psaroudakis@mathematik.uni-stuttgart.de}


\thanks{The first named author is supported by the National Natural Science Foundation of China (Grant No. 11771272). The third named author is supported by Deutsche Forschungsgemeinschaft (DFG, grant KO $1281/14-1$)}


\subjclass[2010]{}

\begin{abstract} 
Recollements of abelian categories are used as a basis of a homological
and recursive approach to quasi-hereditary algebras. 
This yields a homological proof of Dlab and Ringel's characterisation of 
idempotent ideals occuring in heredity chains, which in turn characterises 
quasi-hereditary algebras recursively. Further applications are given to
hereditary algebras and to Morita context
rings.
\end{abstract}

\maketitle


\section{Introduction}

Quasi-hereditary algebras are abundant in representation theory and its
applications to Lie theory and geometry. Examples include hereditary
algebras, Auslander algebras and generally algebras of global dimension two, 
Schur algebras of reductive algebraic groups and other algebras arising from
highest weight categories, endomorphism algebras of projective generators
in categories filtered by standard or exceptional objects, and so on. 
Customary definitions of quasi-hereditary algebras proceed inductively by
first defining what is called a heredity ideal $AeA$ in an algebra $A$ (with
an idempotent element $e=e^2$) and then considering $A/AeA$ and a heredity
ideal therein. The (finite) induction then produces a chain 
$0 \subset Ae_nA \subset Ae_{n-1}A \subset \dots \subset Ae_1A \subset Ae_0A
=A$ with subquotients being heredity ideals in the respective quotient 
algebras. Equivalently, one may define standard modules as being relative
projective over the respective quotient algebra, with $\Delta(n)=Ae_n$,
$\Delta(n-1) = Ae_{n-1}/Ae_n$, and so on. Starting with (semi)simple algebras,
all quasi-hereditary algebras can be constructed using a generalisation
of Hochschild cocycles (see Parshall and Scott's 'not so trivial 
extensions' in \cite{Moosonee}).

Another construction of all quasi-hereditary algebras, recursive in nature
and not using cocycles, has been given by Dlab and Ringel \cite{DR1}, who
were motivated by constructions for perverse sheaves (that are closely
related to quasi-hereditary algebras).
Using ring theoretical methods, Dlab and Ringel 
gave a characterisation of a given
algebra $A$ being quasi-hereditary and a given idempotent ideal $AeA$
occuring somewhere in a heredity chain of $A$, in terms of both $eAe$ and
$A/AeA$ being quasi-hereditary (which is not sufficient) and additional
conditions. 

In the background of all the definitions, characterisations and properties
of quasi-hereditary algebras
are six functors that are the algebraic analogues of
Grothen-dieck's six functors and that form a recollement of abelian categories
relating the module categories of $A$ and of $eAe$ and $A/AeA$. The aim of
this article is to take such recollements and the occuring functors as
basic ingredients for redeveloping the theory of quasi-hereditary algebras,
replacing ring theoretical by homological tools and the inductive approach
(starting with heredity ideals) by a recursive characterisation (starting with
any ideal in a heredity chain) in Theorem~\ref{quai}, which is proved by a direct and homological approach; the main result of \cite{DR1} then 
follows quickly. Another approach via recollements of abelian categories, has been considered by Krause \cite{Krause}; this approach concentrates on heredity ideals.

On the way, various basic properties of quasi-hereditary algebras are given new proofs. Feasibility of the new approach is demonstrated further by
also giving a homological proof that hereditary algebras are quasi-hereditary with any ordering (a result due to Dlab and Ringel) and by adding a class of Morita context rings to the known classes of examples of 
quasi-hereditary algebras. In addition, our approach provides a solution to the problem when the middle term in a recollement of module categories (over semiprimary rings) is hereditary.

\section{Quasi-hereditary algebras and ideals in heredity
chains}

Let $A$ be a semiprimary ring. Let $X$ be a poset and assume that $\{S(x) \ | \ x\in X\}$ is a complete set of pairwise non-isomorphic simple $A$-modules. Semiprimary rings are perfect \cite{Bass}, hence every module has a projective cover. We write $P(x)$ for the projective cover of the simple $A$-module $S(x)$.

Let $N := \{ N_1, \dots, N_k \}$ be a finite set of $A$-modules. The category
of $A$-modules with $N$-filtration, denoted by $\F(_{A}N)$, is defined to be the full subcategory of ${A}\lMod$, i.e. the category of all left $A$-modules, consisting of $A$-modules $M$ such that there exists a filtration 
$0=M_{l+1}\subseteq M_{l}\subseteq \cdots \subseteq M_0=M$ where each quotient 
$M_j/M_{j+1}$ belongs to $\Add{N_i}$ for some $N_i$ ($i$ depending on  $j$)
in $N = \{ N_1, \dots, N_k \}$. Note that the filtration has to be finite,
while the subquotients may be infinite sums in $\Add{N_i}$, which is 
the full subcategory of ${A}\lMod$ consisting of all modules which are summands of a direct sum of $N_i$. Objects in $\F(_{A}N)$ are said to be {\em filtered} 
by $N_1, \dots, N_k$.

Recall from \cite{CPS1} that the ring $A$ is called {\em quasi-hereditary} 
with respect to the poset $X$ if for each $x\in X$, there is a quotient module 
$\Delta(x)$ of $P(x)$, called a {\em standard module}, satisfying the 
following two conditions$\colon$
\begin{enumerate}
\item the kernel of the canonical epimorphism $P(x)\lxr \Delta(x)$ is filtered 
by $\Delta(z)$ with $z>x$, and
\item the kernel of the canonical epimorphism $\Delta(x)\lxr S(x)$ is filtered 
by $S(y)$ with $y<x$.
\end{enumerate}

Recollements  of triangulated or abelian categories were introduced by Beilinson, Bernstein and Deligne in \cite{BBD}. 
A {\em recollement} between abelian categories (see, for instance,
\cite{FP, Kuhn}) $\A,\B$ and $\C$ is a diagram of the form
\[
\xymatrix@C=0.5cm{
\A \ar[rrr]^{\mi} &&&  \B \ar[rrr]^{\me} \ar @/_1.5pc/[lll]_{\mq}  \ar
 @/^1.5pc/[lll]^{{\map}} &&& \C
\ar @/_1.5pc/[lll]_{\ml} \ar
 @/^1.5pc/[lll]^{\mr}
 }
\]
henceforth denoted by $(\A,\B,\C)$, satisfying the following
conditions$\colon$
\begin{enumerate}
\item $(\ml,\me,\mr)$ is an adjoint triple.
\item $(\mq,\mi,\map)$ is an adjoint triple.
\item The functors $\mi$, $\ml$, and $\mr$ are fully faithful.
\item $\Image \mi = \Ker \me$.
\end{enumerate}

For properties of recollements of abelian categories we refer to \cite{FP, Psaroud}. We are interested in recollement with all terms being module categories. Let  $R$ be a ring and let $e$ be an idempotent element of $R$. Then there is a recollement of module categories
\begin{equation}
\label{recolmodulecat}
\xymatrix@C=0.5cm{
{R/ReR}\lMod\ar[rrr]^{\inc} &&& {R}\lMod \ar[rrr]^{eR\otimes_{R}-} \ar @/_1.5pc/[lll]_{R/ReR\otimes_{R}-}  \ar @/^1.5pc/[lll]^{\Hom_{R}(R/ReR,-)} &&& {eRe}\lMod
\ar @/_1.5pc/[lll]_{Re\otimes_{eRe}-} \ar
 @/^1.5pc/[lll]^{\Hom_{eRe}(eR,-)}
 }
\end{equation}
By \cite{PsaroudVitoria}, any recollement of module categories is equivalent, in an appropriate sense, to one induced by
an idempotent element. Thus, the recollement $(\ref{recolmodulecat})$ can be considered as the general recollement situation of ${R}\lMod$. 

If $S$ is a simple $A$-module (resp. simple $eAe$-module), then we denote 
by $P(S)$ (resp. $P_e(S)$) the projective cover of the simple module $S$. 

Now, the main result can be stated and proved; Dlab and Ringel's original
result will follow as Corollary \ref{DRoriginal}. 

\begin{thm}
\label{quai}
Let $A$ be a semiprimary ring and $e$ be an idempotent of $A$. The following statements are equivalent$\colon$
\begin{enumerate}
\item The ring $A$ is quasi-hereditary and there exists a heredity chain such that $AeA$ is contained.

\item There is a recollement of module categories of the form $(\ref{recolmodulecat})$ such that the following conditions hold$\colon$
\begin{enumerate}
\item $A/AeA$ and $eAe$ are quasi-hereditary rings;

\item The counit map $Ae\otimes_{eAe}eA(1-e)\lxr A(1-e)$ is a monomorphism;

\item $eA\in \F(_{eAe}\Delta)$;

\item $\Tor_1^{eAe}(Ae, {_{eAe}\Delta})=0$.
\end{enumerate}
\end{enumerate}
\begin{proof} Let $S_1, S_2, \ldots, S_m, S_{m+1}, \ldots, S_n$ be a full set of non-isomorphic simple $A$-modules. Note that the poset $X$ is now the set $\{1,2,\ldots, n\}$. Indices are chosen such that $eS_i=0$ for all $1\leq i\leq m$ and $eS_i\neq 0$ for all $m+1\leq i\leq n$. Then $\{S_1, \ldots, S_m\}$ are the simple $A/AeA$-modules and $\{eS_{m+1}, \ldots, eS_n \}$ are the simple $eAe$-modules and $A/AeA\otimes_{A}S_{i}=0$ for all $m+1\leq i\leq n$. Moreover, there is an epimorphism $Ae\otimes_{eAe}eS_i\lxr S_i$ for any $m+1\leq i\leq n$. Furthermore, considering the exact sequence $0\lxr AeA\lxr A\lxr A/AeA\lxr 0$ of right $A$-modules and  applying $-\otimes_{A}S_i$, we have that $A/AeA\otimes_{A}S_i=S_i$ and $\Tor^{A}_{1}(A/AeA, S_i)=0$ since $eS_i=0$ for all $1\leq i\leq m$.

Step $0$. If $eS_i \neq 0$, then $eS_i=\Hom_A(Ae,S_i) \neq 0$ and thus $S_i$ is a quotient of $Ae$. Hence there exists a primitive idempotent $e_i$ with $e \cdot e_i = e_i = e_i \cdot e$ such that $Ae_i$ is a projective cover of $S_i$, i.e. isomorphic to $P(S_i)$. Then $eAe_i$ is a projective cover of $eS_i$, thus isomorphic to $P_e(eS_i)$. We use this step later in the proof. 

\vskip 3pt

(i) $\Longrightarrow$ (ii)$\colon$ Let $\Delta(1),\ldots, \Delta(n)$ be the standard $A$-modules up to isomorphism. The proof is divided into seven steps. 

\vskip 3pt

Step $1$. We show that $e\Delta(i)=0$ and $\Tor_1^A(A/AeA, \Delta(i))=0$ for all $1\leq i\leq m$. Since $A$ is quasi-hereditary, there is an exact sequence
\begin{equation}
\label{shortexactseqDeltas}
\xymatrix{
 0 \ar[r]^{} & \Ker{f_i} \ar[r]^{} &   \Delta(i)  \ar[r]^{ \ f_i} & S_i 
 \ar[r] & 0 }
\end{equation}
such that $\Ker f_i$ is filtered by $S_j$ with $1\leq j<i$. Applying the exact functor $eA\otimes_{A}-$ to the filtration of $\Ker f_i$, it follows that $e\Ker f_i=0$. This implies that $e\Delta(i)=0$. 

Consider now the exact sequence $0\lxr AeA\lxr A\lxr A/AeA\lxr 0$ of right $A$-modules. Applying the functor $-\otimes_{A}\Delta(i)$, we get the exact sequence$\colon$ 
\[
\xymatrix@C=0.4cm{
0\ar[r] & \Tor_1^{A}(A/AeA, \Delta(i)) \ar[r] & AeA\otimes_{A}\Delta(i) \ar[r]^{} & A\otimes_{A}\Delta(i) \ar[r]^{}  & A/AeA\otimes_{A}\Delta(i) \ar[r]^{} & 0 }
\]
Since $e\Delta(i)=0$, it follows that $\Tor_1^A(A/AeA, \Delta(i))=0$ for all $1\leq i\leq m$.

\vskip 3pt

Step $2$. We show that $Ae\otimes_{eAe}eP(S_i)\simeq P(S_i)$ and $A/AeA\otimes_{A}\Delta(i)=0$ for all $m+1\leq i\leq n$. By Step $0$, $P_e(eS_i)=eAe_i=eP(S_i)$. Thus, we get an isomorphism $Ae\otimes_{eAe}eAe_i\simeq Ae_i$ showing the first claim. 
For the second one, $A/AeA\otimes_{A}P(S_i)\cong A/AeA\otimes_{A}Ae\otimes_{eAe}eP(S_i)=0$ since $A/AeA\otimes_{A}Ae=0$. Since $A$ is quasi-hereditary, there is an epimorphism $P(S_i)\lxr \Delta(i)$ and therefore $A/AeA\otimes_{A}\Delta(i)=0$.

\vskip 3pt

Step $3$. We show that $eAe$ is a quasi-hereditary ring with standard modules $\{e\Delta({m+1}), \ldots, e\Delta({n})\}$. For every $m+1\leq i\leq n$ there is an exact sequence of the form $(\ref{shortexactseqDeltas})$ such that $\Ker{f_i}$ is filtered by $S_j$ with $1\leq j<i$. Applying the exact functor $eA\otimes_A-$, we obtain the exact sequence $0\lxr e\Ker{f_i}\lxr e\Delta(i)\lxr eS_i\lxr 0$ and $e\Ker{f_i}$ is filtered by $eS_j$ for $m+1\leq j<i$. Also, for all $m+1\leq i\leq n$ there is an exact sequence 
\begin{equation}
\label{shortexactseqsimple}
\xymatrix{
 0 \ar[r]^{} & \Ker{g_i} \ar[r]^{} &   P(S_i)  \ar[r]^{ \ g_i} & \Delta(i) 
 \ar[r] & 0 }
\end{equation}
such that $\Ker g_i$ is filtered by $\Delta(j)$ with $i<j\leq n$.   
Applying the exact functor $e(-)$ we get the exact sequence $(*)\colon 0\lxr e\Ker{g_i}\lxr eP(S_i)\lxr e\Delta(i)\lxr 0$. Note that by Step $0$ the module $eP(S_i) \simeq P_e(eS_i)$ is projective. Clearly, the module $e\Ker{g_i}$ is filtered by $e\Delta(j)$ with $i<j\leq n$. Hence, $eAe$ is a quasi-hereditary ring with standard $eAe$-modules $\{e\Delta(m+1), \ldots, e\Delta(n)\}$. 

Moreover, since the modules $e\Ker{g_i}$ and $e\Delta(i)$ belong to $\F(_{eAe}\Delta)$, the module $eP(S_i)$ lies in $\F(_{eAe}\Delta)$ for all $1\leq i\leq n$ and thus condition (c) holds.

\vskip 3pt

Step $4$. We show that $Ae\otimes_{eAe}e\Delta(i)\simeq \Delta(i)$ and $\Tor^{eAe}_{1}(Ae, e\Delta(i))=0$ for all $m+1\leq i\leq n$  (condition (d)). Consider the canonical morphism $\mu_{\Delta(i)}\colon Ae\otimes_{eAe}e\Delta(i)\lxr \Delta(i)$. The cokernel of $\mu_{\Delta(i)}$ is $A/AeA\otimes_A\Delta(i)$ which is zero by Step $2$. Hence, the map $\mu_{\Delta(i)}$ is an epimorphism for any $m+1\leq i\leq n$. Since $A$ is quasi-hereditary, 
for all $m+1\leq i\leq n$ there is an exact sequence of the form $(\ref{shortexactseqsimple})$ such that $\Ker{g_i}$ is filtered by $\Delta(j)$ with $i<j\leq n$. We claim that the map $\mu_{\Ker{g_i}}\colon Ae\otimes_{eAe}e\Ker{g_i}\lxr \Ker{g_i}$ is an epimorphism. Indeed, let $0\subseteq M_1\subseteq M_2\subseteq\cdots\subseteq M_{n-1}\subseteq  \Ker{g_i}$ be the filtration by $\Delta(j)$ with $i<j\leq n$. Then there are exact sequences
\[
\xymatrix@C=0.5cm{
0 \ar[r]^{} & M_1 \ar[r]^{} & M_2 \ar[r] & M_2/M_1 \ar[r] & 0}, 
\xymatrix@C=0.5cm{0 \ar[r]^{} & M_2 \ar[r]^{} & M_3 \ar[r] & M_3/M_2 \ar[r] & 0}, 
\]
\[           
\cdots, \xymatrix@C=0.5cm{0 \ar[r]^{} & M_{n-1} \ar[r]^{} & \Ker{g_i} \ar[r] & \Ker{g_i}/M_{n-1} \ar[r] & 0
}
\]
such that $M_1$, $M_{t}/M_{t-1}$ for $t=2,\ldots, n-1$, and $\Ker g_i/M_{n-1}$ belong to the set $\{\Delta(i+1), \Delta(i+2), \ldots, \Delta(n)\}$.
Applying $Ae\otimes_{eAe}eA\otimes_{A}-$ to the first exact sequence yields the following exact commutative diagram
\begin{equation}
\label{comdiagram}
\xymatrix@C=0.5cm{
& Ae\otimes_{eAe}eM_1 \ar[r]^{} \ar[d]^{\mu_{M_1}} & Ae\otimes_{eAe}eM_2 \ar[r]^{} \ar[d]^{\mu_{M_2}}  & Ae\otimes_{eAe}e(M_2/M_1)\ar[d]^{\mu_{M_2/M_1}}  \ar[r]^{} & 0 \\
0\ar[r]^{} & M_1 \ar[r]^{} & M_2 \ar[r]^{} & M_2/M_1\ar[r]^{} & 0}
\end{equation}
By diagram chase, the map $\mu_{M_2}$ is an epimorphism. Continuing inductively, with respect to the above exact sequences of the filtration of $\Ker{g_i}$, we obtain that $\mu_{\Ker{g_i}}$ is an epimorphism. Consider now the exact commutative diagram
\[
\xymatrix@C=0.25cm{
 0\ar[r]^{} & \Tor^{eAe}_{1}(Ae, e\Delta(i))\ar[r]^{} & Ae\otimes_{eAe}e\Ker g_i \ar[r]^{} \ar[d]^{} & Ae\otimes_{eAe}eP(S_i) \ar[r]^{} \ar[d]^{}  & Ae\otimes_{eAe}e\Delta(i) \ar[d]^{}  \ar[r]^{} & 0 \\
& 0\ar[r]^{} & \Ker g_i \ar[r]^{} & P(S_i) \ar[r]^{} & \Delta(i)\ar[r]^{} & 0}
\]
Step $2$ provides an isomorphism $Ae\otimes_{eAe}eP(S_i)\simeq P(S_i)$ for all $m+1\leq 1\leq n$. Then, by Snake Lemma and since the map $\mu_{\Ker{g_i}}$ is an epimorphism, we have $Ae\otimes_{eAe}e\Delta(i)\simeq \Delta(i)$. Since $\Ker g_i$ is filtered by $\Delta(j)$ with $i<j\leq n$, we also get that $Ae\otimes_{eAe}e\Ker g_i\simeq \Ker g_i$ for all $m+1\leq i\leq n$. This implies that $\Tor^{eAe}_{1}(Ae, e\Delta(i))=0$ for all  $m+1\leq i\leq n$.

\vskip 3pt

Step $5$. We show that the map $\mu_{P(S_i)}\colon Ae\otimes_{eAe}eP(S_i)\lxr P(S_i)$ is a monomorphism for all $1\leq i\leq m$. Consider the short exact sequence $(\ref{shortexactseqsimple})$. Since $e\Delta(i)=0$ for all $1\leq i\leq m$ by Step $1$, there is the following exact commutative diagram$\colon$
\begin{equation}
\label{comdiagrammono}
\xymatrix{
& 0\ar[r] & Ae\otimes_{eAe}e\Ker{g_i} \ar[r]^{\cong} \ar[d]^{\mu_{\Ker{g_i}}} & Ae\otimes_{eAe}eP(S_i) \ar[r]^{} \ar[d]^{\mu_{P(S_i)}}  & 0   &  \\
& 0\ar[r]^{} & \Ker g_i\ar[r]^{} & P(S_i) \ar[r]^{g_i} & \Delta(i) \ar[r] & 0 }
\end{equation}
We claim that the map $\mu_{\Ker{g_i}}$ is a monomorphism. Consider the filtration of $\ker{g_i}$, and in particular, the exact sequence $0\lxr M_1\lxr M_2\lxr M_2/M_1\lxr 0$ (as in the proof of Step $4$) where $M_1$ and $M_2/M_1$ lie in the set $\{\Delta(i+1), \ldots, \Delta(n)\}$. 
Note that $j>i$ and $1\leq i\leq m$. Consider now the diagram $(\ref{comdiagram})$. Step $1$ implies $e\Delta(i)=0$ for all $1\leq i\leq m$ and by Step $4$ we have $Ae\otimes_{eAe}e\Delta(i)\simeq \Delta(i)$ for all $m+1\leq i\leq n$. Clearly, in any of the latter cases the map $\mu_{M_2}$ in $(\ref{comdiagram})$ is a monomorphism. Continuing inductively on the length of the filtration of $\Ker{g_i}$, the map $\mu_{\Ker{g_i}}$ is seen to be a monomorphism. Then from diagram $(\ref{comdiagrammono})$ it follows that the map $\mu_{P(S_i)}$ is a monomorphism for any $1\leq i\leq m$, that is, condition (b) holds.

\vskip 3pt

Step $6$. We show that $\Tor^{A}_{1}(A/AeA, \Delta(i))=0$ for all $m+1\leq i\leq n$. Consider the exact sequence $0\lxr AeA\lxr A\lxr A/AeA\lxr 0$ of right $A$-modules. By Step $2$ we have the following exact sequence$\colon$ 
\[
\xymatrix{
0\ar[r] & \Tor_1^{A}(A/AeA, \Delta(i)) \ar[r] & AeA\otimes_{A}\Delta(i) \ar[r]^{} & A\otimes_{A}\Delta(i) \ar[r]^{}  & 0 }
\]
Consider the following exact commutative diagram$\colon$
\[
\xymatrix{
Ae\otimes_{eAe}eA \ar[rr]^{\mu_{A}} \ar@{->>}[dr]_{\kappa_A} && A \ar[r] & A/AeA \ar[r] & 0 \\
  & AeA \ \ \ \ar@{>->}[ru]_{\lambda_A} & & &
}
\]
Applying the functor $-\otimes_{A}\Delta(i)$, we get the commutative diagram
\[
\xymatrix{
Ae\otimes_{eAe}eA\otimes_A\Delta(i) \ar[rr]^{\mu_{A}\otimes \Delta(i)} \ar@{->>}[dr]_{\kappa_A\otimes \Delta(i)} && A\otimes \Delta(i) \ar[r] & 0 \\
  & AeA\otimes_A\Delta(i) \ar[ru]_{\lambda_A\otimes \Delta(i)} & & 
}
\]
From Step $4$, the map $\mu_{A}\otimes \Delta(i)$ is an isomorphism. This implies that the map $\kappa_A\otimes \Delta(i)$ is an isomorphism and the map $\lambda_A\otimes \Delta(i)$ is an epimorphism. By the commutativity of the above diagram, we get the desired $\Tor$-vanishing.

\vskip 3pt

Step $7$. We show that the ring $A/AeA$ is quasi-hereditary with standard modules $\{A/AeA\otimes_{A}\Delta(1), \ldots, A/AeA\otimes_{A}\Delta(m)\}$. Recall that for all $1\leq i\leq m$ we have $\Tor^{A}_{1}(A/AeA, S_i)=0$ (see the first paragraph of the proof). Applying the functor $A/AeA\otimes_{A}-$ to the short exact sequence $(\ref{shortexactseqDeltas})$, we get the short exact sequence $0\lxr A/AeA\otimes_{A}\Ker f_i\lxr A/AeA\otimes_{A}\Delta(S_{i})\lxr A/AeA\otimes_{A}S_i\lxr 0$ such that $A/AeA\otimes_{A}\Ker f_i$ is filtered by $A/AeA\otimes_A S_j$ ($\cong S_j$) with $1\leq j<i$. Consider now the short exact sequence $(\ref{shortexactseqsimple})$. Recall that in this case $\Ker{g_i}$ is filtered by $\Delta(S_j)$ with $i<j\leq n$. Then, by Step $6$ we obtain the exact sequence 
$0\lxr A/AeA\otimes_{A}\Ker{g_i}\lxr A/AeA\otimes_{A}P(S_i)\lxr A/AeA\otimes_{A}\Delta(i)\lxr 0$ such that $A/AeA\otimes_{A}\Ker{g_i}$ is filtered by $A/AeA\otimes_{A}\Delta(j)$ with $i<j\leq n$. We infer that the ring $A/AeA$ is quasi-hereditary.

\vskip 3pt

(ii) $\Longrightarrow$ (i)$\colon$ 
Since the ring $eAe$ is quasi-hereditary, there exist standard $eAe$-modules
\[
\big\{\Delta(eS_{m+1}), \ldots, \Delta(eS_n)\big\}.
\]
Also, since the ring $A/AeA$ is quasi-hereditary, there are standard $A/AeA$-modules
\[
\big\{\Delta(S_{1}), \ldots, \Delta(S_m)\big\}.
\]
The proof is divided into two steps. 

\vskip 3pt

Step $1$. We show that $\{Ae\otimes_{eAe}\Delta(eS_{m+1}), \ldots, Ae\otimes_{eAe}\Delta(eS_n)\}$ are standard $A$-modules. First, recall that $Ae\otimes_{eAe}eP(S_i)\simeq P(S_i)$ for any $m+1\leq i\leq n$, see the proof of  Step $2$  in (i) $\Longrightarrow$ (ii).
The latter isomorphism together with condition (b)
gives the isomorphism $Ae\otimes_{eAe}eA\simeq AeA$.

Since $eAe$ is quasi-hereditary,  there exists an exact sequence
\begin{equation}
\label{exa3}
\xymatrix{
  0\ar[r]^{} & \Ker\phi_i \ar[r]^{ } & \Delta(eS_i) \ar[r]^{\ \ \ \phi_i \ \  } & eS_i  \ar[r]^{} & 0}
\end{equation}
for all $m+1\leq i\leq n$, such that $\Ker\phi_i$ is filtered by $eS_j$ with $m+1\leq j<i$. Applying $Ae\otimes_{eAe}-$ to (\ref{exa3}) gives an exact sequence
\[
\xymatrix{
  0\ar[r]^{} & \Ker \psi_i \ar[r]^{ } & Ae\otimes_{eAe}\Delta(eS_i) \ar[r]^{\ \ \ \ \ \ \ \ \ \ \ \psi_{i} \ \  } & S_i  \ar[r]^{} & 0}
\]
Moreover, $e\Ker \psi_i\simeq \Ker\phi_i$. 
We claim that $\Ker \psi_i$ is filtered by $S_j$ with $1\leq j<i$. Assume to the contrary that
 $0\subseteq L_1\subseteq L_2\subseteq\cdots\subseteq L_{n-1}\subseteq  \Ker\psi_i$ is a filtration of $\Ker\psi_i$  by $S_j$ for $1\leq j\leq n$. Then there are exact sequences $0\lxr L_1\lxr L_2\lxr L_2/L_1\lxr 0$, $0\lxr L_2\lxr L_3\lxr  L_3/L_2\lxr 0$ and so on, where $L_1$ and all the quotients are simple $A$-modules.
Applying $eA\otimes_{A}-$ to the above sequences,  we obtain that $\Ker\phi_i$ is filtered by $eS_j$ with $m+1\leq j\leq n$, which is a contradiction since $j$ is strictly smaller that $i$. Hence, our claim holds.

On the other hand, in the short exact sequence
\begin{equation}
\label{exa1}
\xymatrix{
  0\ar[r]^{} & \Ker{f_i} \ar[r]^{ } & P_{e}(eS_i) \ar[r]^{\ \ \ f_i \ \  } & \Delta(eS_i)  \ar[r]^{} & 0}
\end{equation}
the first term $\Ker{f_i}$ is filtered by $\Delta(eS_j)$ for $i<j\leq n$. Since $\Tor_1^{eAe}(Ae, \Delta(eS_i))=0$ for all $m+1\leq i\leq n$ by condition (d), applying the functor $Ae\otimes_{eAe}-$ to $(\ref{exa1})$ yields the short exact sequence$\colon$
\[
\xymatrix{
  0\ar[r]^{} & Ae\otimes_{eAe}\Ker{f_i} \ar[r]^{ } & P(S_i) \ar[r]^{ } & Ae\otimes_{eAe}\Delta(eS_{i})  \ar[r]^{} & 0}
\]
such that $Ae\otimes_{eAe}\Ker{f_i}$ is filtered by $Ae\otimes_{eAe}\Delta(eS_{j})$ with $i<j\leq n$.

\vskip 3pt

Step $2$. We prove that $\{\Delta(S_{1}), \ldots, \Delta(S_{m})\}$ are standard $A$-modules. Since $A/AeA$ is quasi-hereditary, for all $1\leq i\leq m$ there is an exact sequence of left $A/AeA$-modules, and thus also of left $A$-modules,
\[
\xymatrix{
  0\ar[r]^{} & \Ker\varphi_i \ar[r]^{ } & \Delta(S_{i}) \ar[r]^{ \ \ \varphi } & S_i  \ar[r]^{} & 0}
\]
such that $\Ker\varphi_i$ is filtered by $S_j$ with $1\leq j<i$.  
Since  $A/AeA\otimes_{A}P(S_i)=A/AeA\otimes_{A}Ae_i=(A/AeA)e_i$,
it follows that $A/AeA\otimes_{A}P(S_i)$ is the projective cover of $S_i$ as an $A/AeA$-module. Consider the epimorphism $h_i\colon A/AeA\otimes_{A}P(S_i)\lxr \Delta(S_i)$ where $\ker h_i$ is filtered by $\Delta(S_j)$ for $i< j\leq m$. By assumption (b), for all $1\leq i\leq m$ there is the following short exact sequence of left $A$-modules
\[
\xymatrix{
  0\ar[r]^{} & Ae\otimes_{eAe}eP(S_i) \ar[rr]^{\ \ \ \ \ \mu_{P(S_i)}} && P(S_i) \ar[rr]^{\lambda_{P(S_i)} \ \ \ \ \ \ \ \ } && A/AeA\otimes_{A}P(S_i)  \ar[r]^{} & 0}
\]
Define the composition $g_i:=h_i\circ\lambda_{P(S_i)}$ and consider the short exact sequence
\[
\xymatrix{
  0\ar[r]^{} & \Ker{g_i} \ar[r]^{} & P(S_i) \ar[r]^{\ \ \ \ \ \ \  g_i\ \ \ \  \ \ \ } & \Delta(S_i)  \ar[r]^{} & 0}
\]
We claim that the first term $\Ker{g_i}$ is filtered by $\Delta(S_j)$ for $i< j\leq n$.  
Applying the Snake Lemma to the commutative diagram
\[
\xymatrix{
 & Ae\otimes_{eAe}eP(S_i) \ar@{=}[d]^{} \ar[rr]^{} && \Ker g_i \ar[rr]^{} \ar@{^{(}->}[d]^{\inc} && \Ker h_i  \ar@{^{(}->}[d]^{\inc}  \ar[r]^{} & 0 
 \\
 0\ar[r]^{} & Ae\otimes_{eAe}eP(S_i) \ar[rr]^{\ \ \ \ \ \mu_{P(S_i)}} && P(S_i) \ar[rr]^{\lambda_{P(S_i)} \ \ \ \ \ \ } \ar@{->>}[d]^{g_i} && A/AeA\otimes_{A}P(S_i) \ar@{->>}[d]^{h_i}  \ar[r]^{} & 0 \\
   & && \Delta(S_i) \ar@{=}[rr] && \Delta(S_i) & }
\]
provides us with the short exact sequence
\begin{equation}
\label{lastses}
\xymatrix{
  0\ar[r]^{} & Ae\otimes_{eAe}eP(S_i) \ar[r]^{} & \Ker{g_i} \ar[r]^{ } & \Ker{h_i}  \ar[r]^{} & 0}
\end{equation}
By assumption (c) and (d), for all $1\leq i\leq m$, the $eAe$-module $eP(S_i)$ is filtered by $\Delta(eS_j)$ for $m+1\leq j\leq n$ and $\Tor^{eAe}_1(Ae,\Delta(eS_j))=0$. Therefore, $Ae\otimes_{eAe}eP(S_i)$ is filtered by $Ae\otimes_{eAe}\Delta(eS_j)$ with $m+1\leq j\leq n$. Moreover, the module $\Ker{h_i}$ is filtered by $\Delta(S_j)$ for $i< j\leq m$. From 
$(\ref{lastses})$ it follows that $\Ker{g_i}$ is filtered by $Ae\otimes_{eAe}\Delta(eS_j)$ and $\Delta(S_j)$ for $i<j\leq n$.

\vskip 3pt

By Step $1$ and Step $2$, the ring $A$ is quasi-hereditary.
\end{proof}
\end{thm}

\begin{rem}
\label{remoppositering}
Dlab and Ringel have shown that $A$ is a quasi-hereditary ring if and only if the opposite ring $A^{\op}$ is quasi-hereditary \cite[Statement 9]{DR2}. The proof proceeds inductively and is based on the fact that for a heredity ideal $AeA$ in a ring, multiplication in $A$ provides an isomorphism $(\dagger) \ Ae \otimes_{eAe} eA \stackrel{\mathsf{mult}}{\lxr} AeA$ of
$A$-bimodules. When $AeA$ is a heredity ideal, then $eAe$ is semisimple, and then
multiplication in $(\dagger)$ is an isomorphism if and only if $AeA$ is
projective as a left $A$-module if and only if it is projective as a right
$A$-module. This is the left-right symmetry needed. The isomorphism
$(\dagger)$ is a special case of a direct consequence of condition (b) in Theorem~\ref{quai}. Using the result of Dlab and Ringel, Theorem~\ref{quai} can be reformulated for the ring $A^{\op}$ in terms of conditions (a)$^{\op}$--(d)$^{\op}$.
\end{rem}

Let $A$ be a quasi-hereditary ring with heredity chain $\mathcal{J}=(J_i)_{0\leq i\leq n}$ and let $M$ be an $A$-module. Then the $\mathcal{J}$-filtration of $M$
is the chain of submodules
\[
0=J_{n+1}M\subseteq J_nM\subseteq \cdots \subseteq J_0M=M.
\]
Dlab and Ringel \cite{DR1} called the $\mathcal{J}$-filtration of $M$ {\em good} if the quotient $J_iM/J_{i+1}M$ is projective as an $A/J_{i+1}$-module for all $0\leq i\leq n$.

\begin{lem}
\label{lemDeltafiltered}
Let $A$ be a quasi-hereditary ring and let $M$ be a left $A$-module. The following statements are equivalent.
\begin{enumerate}
\item The $\mathcal{J}$-filtration of $M$ is good.

\item $M\in \F(_A\Delta)$.

\end{enumerate}
\begin{proof}
The heredity ideal $J_n$ is generated by a primitive idempotent $e_n$, hence
$J_n = A e_nA$. Therefore, $J_nM = Ae_nM$ is the trace of $Ae_n$ in $M$. It is
projective if and only if it is in $\Add(Ae_n)$, which equals $\Add(\Delta(n))$
because of $\Delta(n) \simeq Ae_n$. As $\Delta(j)$ for $j \neq n$ has no
composition factor $S(j)$ and hence
$\Hom_A(Ae_n,\Delta(j)) = e_n \Delta(j) = 0$ for all
$j \neq n$, the bottom part of a $\Delta$-filtration of $M$, if there is one,
must coincide with $Ae_nM$, which then must be projective. Continuing by
induction, the claimed equivalence follows.
\end{proof}
\end{lem}

As a consequence of Theorem~\ref{quai}, Remark~\ref{remoppositering} and Lemma~\ref{lemDeltafiltered} we obtain the main result of Dlab and Ringel \cite{DR1}: 

\begin{cor}{\textnormal{(\cite[Theorem 1]{DR1})}} \label{DRoriginal}
Let $A$ be a semiprimary ring 
and $e$ an idempotent element of $A$. The following statements are equivalent$\colon$
\begin{enumerate}
\item There is a heredity chain for $A$ containing $AeA$.

\item The rings $A/AeA$ and $eAe$ are quasi-hereditary, the multiplication map
\[
\xymatrix{
  Ae\otimes_{eAe}eA\ar[r]^{} & AeA}
\]
is bijective, and there is a heredity chain $\mathscr{I}$ of $eAe$ such that the $\mathscr{I}$-filtrations of $Ae_{eAe}$ and ${_{eAe}eA}$ are good.

\item The rings $A/AeA$ and $eAe$ are quasi-hereditary, the multiplication map
\[
\xymatrix{
  (1-e)Ae\otimes_{eAe}eA(1-e)\ar[r]^{} & (1-e)AeA(1-e)}
\]
is bijective, and there is a heredity chain $\mathscr{I}$ of $eAe$ such that the $\mathscr{I}$-filtrations of $(1-e)Ae_{eAe}$ and ${_{eAe}eA(1-e)}$ are good.
\end{enumerate}
\end{cor}

\section{Further uses of the homological approach}

In this section we provide several applications of Theorem~\ref{quai}. We start by showing that quasi-hereditary algebras with any ordering coincide with the class of hereditary algebras. This is a classical result due to Dlab and Ringel. Using Theorem~\ref{quai} we give a new proof which simultaneously provides an answer to the following problem on hereditary rings in a recollement situation.

Let $(\A,\B,\C)$ be a recollement of abelian categories. By \cite[Theorem 4.8]{Psaroud:homolrecol}, if $\B$ is hereditary, i.e. $\gld\B\leq 1$, then $\A$ and $\C$ are also hereditary. The converse is wrong, when just one recollement is used. There is, however, a converse in terms of a set of recollements related with heredity chains in semiprimary rings. To state this result, some notation
has to be fixed. Let $A$ be a semiprimary ring and let $X$ be the set
of isomorphism classes of simple $A$ modules. Suppose $X=X_1\sqcup X_2$ is a disjoint union of two non-empty subsets. Let $e_{X_1}$ be an idempotent
such that $Ae_{X_1}$ is a direct sum of projective covers of simple
modules representing all classes in $X_1$, and $e_{X_2}$ similarly.

\begin{cor}
{\textnormal{(\cite[Theorem 1]{DR2})}}
\label{gld1}
Let $A$ be a semiprimary ring. The following statements are equivalent$\colon$ 
\begin{enumerate}
\item $A$ is a hereditary ring. 

\item $A$ is a quasi-hereditary ring with any ordering.

\item For all partitions of $X$ into $X_1\sqcup X_2$, the ring $A$ has a heredity chain such that $Ae_{X_2}A$ is contained and $A/Ae_{X_2}A$, $e_{X_2}Ae_{X_2}$ are hereditary.
\end{enumerate}
\begin{proof} 
(i) $\Longrightarrow$ (ii)$\colon$ Suppose that $A$ is hereditary and let $e$ be a primitive idempotent of $A$. Associated with any idempotent element $e$ we always have a recollement of module categories, see diagram~$(\ref{recolmodulecat})$. Since $\gld{A}\leq 1$, it follows that  $AeA$ is a  projective left $A$-module. Moreover, let $f\colon Ae\lxr Ae$ be a non-zero $A$-morphism. 
We claim that $f$ is an isomorphism. Indeed, if $f$ is an epimorphism then
it is an isomorphism since $Ae$ is indecomposable. Suppose $f$ is not
surjective. If $f$ is not a monomorphism, then $\Ker{f}$ is projective since 
$A$ is hereditary and therefore $\pd\Coker{f}=2$, contradicting the fact that 
$A$ is hereditary. Thus, the morphism $f$ must be a monomorphism and its image
must be contained in $\rad(Ae)$. Hence, $f$ restricts to injective maps
$\rad^j(Ae)\lxr\rad^j(\Image(f)) \subset \rad^{j+1}(Ae)$ for all $j$, a contradiction to $A$ being semiprimary and thus having finite radical length.
This implies $e(\rad{A})e=0$ and therefore the ideal $AeA$ is heredity. Moreover, since $e(\rad{A})e=\rad eAe$, the algebra $eAe$ is semisimple. Thus, the algebra $eAe$ is quasi-hereditary and the conditions (c) and (d) of Theorem~\ref{quai} hold. On the other hand, the projectivity of $AeA$ implies that it is a stratifying ideal, i.e. $Ae\otimes_{eAe}eA\cong AeA$ and $\Tor^i_{eAe}(Ae,eA)=0$ for all $i>0$, thus condition (b) of Theorem~\ref{quai} holds. Since $AeA$ is a stratifying ideal, $\gld{A/AeA}\leq \gld{A}\leq 1$. By induction on the number of simple modules and since $A/AeA$ is hereditary, $A/AeA$ is quasi-hereditary. By Theorem~\ref{quai}, the algebra $A$ is quasi-hereditary with any ordering since $e$ was an arbitrary idempotent element of $A$.

(ii) $\Longrightarrow$ (i)$\colon$ Suppose that $A$ is a quasi-hereditary algebra with any ordering. We proceed by induction on the number of simple modules. Induction starts with a local quasi-hereditary ring $A$. Then $A$ equals a heredity ideal $AeA$, for some idempotent $e$ that must be
equivalent to the unit of $A$. Then $eAe$, which is semisimple, is Morita
equivalent to $A$. Hence $A$ is simple. Assume that we have two non-trivial idempotents $e_1$ and $e_2$ (i.e. $e_2=1-e_1$). Then we have the recollement of module categories $({A/AeA}\lMod,{A}\lMod{A},{eAe}\lMod)$ and we iterate like this.

We continue now by showing that $A$ is hereditary. Let $S$ be a simple $A$-module which is annihilated by a primitive idempotent element $e$ of $A$. 
Since $A$ is a quasi-hereditary algebra with any ordering, it follows that $AeA$ is a heredity ideal and $A/AeA$ is a quasi-hereditary algebra with any ordering. By induction hypothesis the algebra $A/AeA$ is hereditary and therefore we have $\pd_{A/AeA}S\leq 1$. Since $AeA$ is a projective left $A$-module, it follows that $\pd_{A}A/AeA=1$. This implies that $\pd_{A}S\leq 2$. We claim that $\pd_{A}S=2$ is not the case. Since $eS=0$ and $AeA$ is a projective left $A$-module, applying the functor $-\otimes_{A}S$ to the exact sequence $0\lxr AeA\lxr A\lxr A/AeA\lxr 0$ of right $A$-modules, we get that $\Tor^{A}_{n}(A/AeA, S)=0$ for any $n\geq 1$. Let
\[
\xymatrix{
0 \ar[r] & P_2 \ar[r]^{f_2}  & P_1 \ar[r]^{f_1} & P_0 \ar[r]^{f_0} & S \ar[r] & 0  }
\]
be a minimal projective resolution of $S$. Since $\Tor^{A}_{1}(A/AeA, \Ker f_0)=0$, applying the functor $A/AeA\otimes_{A}-$ we obtain the following exact sequence
\[
\xymatrix{
0 \ar[r] & A/AeA\otimes_AP_2 \ar[r]^{\iden\otimes f_2 \ }  & A/AeA\otimes_AP_1 \ar[r]^{\iden\otimes f_1 \ } & A/AeA\otimes_AP_0 \ar[r]^{ \ \ \ \ \ \ \ \iden\otimes f_0} & S \ar[r] & 0  }
\]
where $A/AeA\otimes_AP_i$ are projective left $A/AeA$-modules.
Since $\pd_{A/AeA}S\leq 1$ it follows that either $A/AeA\otimes_{A}P_2=0$ or $\iden\otimes_{A}f_2$ is a non-zero split monomorphism.

\vskip 3pt

First case: $A/AeA\otimes_{A}P_2=0$. Since $AeA\otimes_{A}P_2\simeq P_2$ and $AeA$ is a projective left $A$-module, we get a split exact sequence $0\lxr P_2\lxr AeA\otimes_{A}P_1 \lxr AeA\otimes_{A}P_0 \lxr 0$. Note that $AeA\otimes_AP_1$ is a direct summand of $P_1$, and thus from the splitting we get that $P_2$ is a direct summand of $P_1$. However, this contradicts the minimality of the projective resolution of $S$.

\vskip 3pt

Second case: $\iden\otimes_{A}f_2$ is a non-zero split monomorphism. Let $\iden\otimes_A h$ be the inverse and denote by $\pi_2\colon P_2\lxr A/AeA\otimes_{A}P_2$ the canonical epimorphism. Since $(\iden\otimes_{A}hf_2)\pi_2=\pi_2 (hf_2)$, it follows that
$\pi_2=\pi_2 (hf_2)$. Consider now the following diagram with exact rows$\colon$
\[
\xymatrix{
0 \ar[r] & AeA\otimes_{A}P_2 \ar[r]^{}  & P_2 \ar[r]^{\ \ \ \ \ \ \ \ \ } \ar[d]^{\iden_{P_2}-hf_2} & A/AeA\otimes_{A}P_2 \ar[r]^{} & 0    \\
0 \ar[r] &  AeA\otimes_{A}P_2 \ar[r]^{\ \ \ \ \iota_2} & P_2\ar[r]^{\pi_2\ \ \ \ \ \ \ \ \ } & A/AeA\otimes_{A}P_2 \ar[r]^{} & 0 ,}
\]
Then the map $\iden_{P_2}-hf_2$ factors through $\iota_2$, i.e. there is a map $\psi\colon P_2\lxr AeA\otimes_{A}P_2$ such that $\iden_{P_2}=hf_2+\iota_2\psi$. Hence, $P_2$ is a direct summand of $P_1\oplus (AeA\otimes_{A}P_2)$ and this implies that $P_2$ and $P_1$ have common direct summands. This contradicts the minimality of the projective resolution of $S$. Thus $\pd_{A}S\leq 1$, i.e. $A$ is hereditary.

(ii) $\Longleftrightarrow$ (iii)$\colon$ Assume that (ii) holds and let $X=X_1\sqcup X_2$ be a partition of $X$. Then the ring $A$ has a heredity chain such that $Ae_{X_2}A$ is contained and since $A$ is hereditary, it follows from \cite[Theorem 4.8]{Psaroud:homolrecol} that the rings $A/Ae_{X_2}A$ and $e_{X_2}Ae_{X_2}$ are hereditary. The implication (iii) $\Longrightarrow$ (ii) is clear.
\end{proof}
\end{cor}

Next we provide a sufficient condition for a class of  Morita context rings to be quasi-hereditary. For more details on Morita context rings,
we refer to \cite{GrP}.

\begin{cor}
\label{Mo}
Let $A$ be a finite dimensional $k$-algebra over a field $k$, and $e$ and $f$ two idempotent elements of $A$ such that $fAe=0$. Let
$N:=Ae\otimes_{k}fA$ and $\Lambda_{(0,0)}:=\bigl(\begin{smallmatrix}
A & N \\
N & A
\end{smallmatrix}\bigr)$. If $A$ is a quasi-hereditary algebra, then the Morita context ring $\Lambda_{(0,0)}$ is a quasi-hereditary algebra.
\begin{proof} Since $fAe=0$, it follows that $N\otimes_{A}N=0$. Then by \cite[Example 4.16]{GaP}, $\Lambda_{(0,0)}$ is a Morita context ring, whose addition is componentwise, and multiplication is given as follows$\colon$
\[
         \begin{pmatrix}
           a & n \\
           m & b \\
         \end{pmatrix}
       \cdot
         \begin{pmatrix}
           a' & n' \\
           m' & b' \\
         \end{pmatrix}=
         \begin{pmatrix}
           aa' & an'+nb' \\
           ma'+bm' & bb' \\
         \end{pmatrix}
\]
The objects of $\smod\Lambda_{(0,0)}$ are given by tuples $(X,Y,f,g)$, where $X\in \smod{A}$, $Y\in \smod{A}$, $f\colon N\otimes_AX\lxr Y$ and $g\colon N\otimes_AY\lxr X$. The compatibility conditions that objects over a Morita context ring should satisfy are trivial since $N\otimes_{A}N=0$, see \cite{GrP}. Furthermore, from \cite[Proposition 2.4]{GrP} there is a recollement 
\[
\xymatrix@C=0.5cm{
{A}\lsmod \ar[rrr]^{} &&& \Lambda_{(0,0)}\lsmod \ar[rrr]^{\mU_A} \ar
@/_1.5pc/[lll]_{}  \ar
 @/^1.5pc/[lll]^{} &&& {A}\lsmod
\ar @/_1.5pc/[lll]_{\mt_{A}} \ar
 @/^1.5pc/[lll]^{\mh_{A}}
 }
 \]
\noindent  where $\mt_{A}(X)=(X, N\otimes_{A}X, \iden_{N\otimes_{A}X}, 0)$, $\mU_{A}(X, Y, f, g)=X$ and $\mh_{A}(X)=(N\otimes_{A}X, X, 0,\iden_{N\otimes_{A}X})$. From \cite[Proposition 3.1]{GrP} the indecomposable projective $\Lambda_{(0,0)}$-modules are of the form $\mt_A(P)$ and $\mh_A(P)$, where $P$ is an indecomposable projective $A$-module. We use Theorem~\ref{quai} to derive that $\Lambda_{(0,0)}$ is a quasi-hereditary algebra. The recollement of ${\Lambda_{(0,0)}}\lsmod$ induced by the idempotent element $\varepsilon=\bigl(\begin{smallmatrix}
1 & 0\\
0 & 0
\end{smallmatrix}\bigr)$ is precisely the one given above (consider the recollement $(\ref{recolmodulecat})$ for finitely generated modules). Condition (a) of Theorem~\ref{quai} is clearly satisfied since $A$ is quasi-hereditary. To check condition (b), we compute the counit map of the adjunction $(\mt_A,\mU_A)$. In particular, there are morphisms
\[
\xymatrix{
\mt_A\mU_A(\mt_A(P)) \ar[rrr]^{ \ \ \ (\iden_P,\iden_{N\otimes_AP})} &&& \mt_A(P)
}
\]
and
\[
\xymatrix{
\mt_A\mU_A(\mh_A(P)) \ar[rrr]^{ \ \ \ (\iden_{N\otimes_AP},0)} &&& \mh_A(P)
}
\]
where $\mt_A\mU_A(\mh_A(P))=(N\otimes_AP,0,0,0)$. Hence, the counit map in any projective is a monomorphism, so condition (b) holds.

For conditions (c) and (d) of Theorem~\ref{quai}, observe that $\Lambda_{(0,0)}\varepsilon$ is a projective right $\varepsilon\Lambda_{(0,0)}\varepsilon$-module and $\varepsilon\Lambda_{(0,0)}$ is a projective left $\varepsilon\Lambda_{(0,0)}\varepsilon$-module since $N$ is a both left and right projective $A$-module. Note that $A\simeq \varepsilon\Lambda_{(0,0)}\varepsilon$.

By Theorem~\ref{quai}, $\Lambda_{(0,0)}$ is quasi-hereditary.
\end{proof}
\end{cor}

Let $A$ now be a finite dimensional quasi-hereditary algebra over a field $k$ and with respect to a poset $X$. The $k$-duals of the standard $A^{\op}$-modules are $A$-modules, which are called costandard. Recall
from \cite{CPS1} that for each $x\in X$, the costandard module $\nabla(x)$ satisfies the following two conditions$\colon$
\begin{enumerate}
\item there is a monomorphism $L(x)\lxr \nabla(x)$ such that the cokernel  is filtered by $L(y)$ with $y<x$;

\item there is a monomorphism $\nabla(x)\lxr I(x)$ such that the cokernel is filtered by $\nabla(z)$ with $z>x$.
\end{enumerate}
We denote by $\F(_{A}\nabla)$ the full subcategory of ${A}\lsmod$ consisting of $A$-modules which have a filtration by costandard $A$-modules.

Ringel \cite{R} introduced the notion of the characteristic tilting module, which is a basic module $T$ such that $\F(_{A}\Delta)\cap \F(_{A}\nabla)=\add{T}$. 
We close this section with the next result, where we investigate the behaviour of the characteristic tilting module along the recollement situation $(\ref{recolmodulecat})$ of Theorem~\ref{quai}. We remark that we consider below a version of $(\ref{recolmodulecat})$ for finitely generated modules.

\begin{cor}
\label{cha}
Let $A$ be a quasi-hereditary algebra such that $AeA$ is contained in a heredity chain of $A$. The following hold.
\begin{enumerate}

\item The functor $Ae\otimes_{eAe}-\colon {eAe}\lsmod\lxr {A}\lsmod$ sends $\F(_{eAe}\Delta)$ to $\F(_{A}\Delta)$.

\item The functor $eA\otimes_{A}-\colon {A}\lsmod\lxr {eAe}\lsmod$ sends $\F(_{A}\Delta)$, resp. $\mathcal{F}({_{A}\nabla})$, to $\F(_{eAe}\Delta)$, resp. $\mathcal{F}({_{eAe}\nabla})$.

\item The inclusion functor $\inc\colon {A/AeA}\lsmod\lxr {A}\lsmod$ sends $\F(_{A/AeA}\Delta)$, resp. $\mathcal{F}({_{A/AeA}\nabla})$, to $\F(_{A}\Delta)$, resp. $\mathcal{F}({_{A}\nabla})$.

\item The functors $eA\otimes_{A}-$ and $\inc$ preserve the characteristic tilting modules.
\end{enumerate}
\begin{proof} (i) This follows immediately using condition (d), i.e. $\Tor_1^{eAe}(Ae, {_{eAe}\Delta})=0$, of Theorem~\ref{quai}.

(ii) First, from the proof of Step $3$ in (i) $\Longrightarrow$ (ii) of Theorem~\ref{quai}, we have that the functor $eA\otimes_{A}-\colon {A}\lsmod\lxr {eAe}\lsmod$ sends $\mathcal{F}(_{A}\Delta)$ to $\mathcal{F}({_{eAe}\Delta})$. We show that the functor $eA\otimes_{A}-$ sends $\mathcal{F}(_{A}\nabla)$ to $\mathcal{F}(_{eAe}\nabla)$. Indeed, since $A$ is  quasi-hereditary, the opposite algebra $A^{\op}$ is also quasi-hereditary. Denote by $\mathsf{D}=\Hom_{k}(-.k)$ the standard $k$-duality and let 
$\nabla(S_1), \ldots, \nabla(S_{n})$ be all costandard $A$-modules. Then $\mathsf{D}(\nabla(S_i))$ is a standard $A^{\op}$-module for each $1\leq i\leq n$. Thus, we get that $\mathsf{D}(\nabla(S_i))e$ is a standard $(eAe)^{\op}$-module. Since $\mathsf{D}(e\nabla(S_i))\cong \mathsf{D}(\nabla(S_i))e$, it follows that $e\nabla(S_i)$ is a costandard $eAe$-module for each $1\leq i\leq n$. 

(iii) By Step $2$ of (ii) $\Longrightarrow$ (i) in Theorem~\ref{quai}, the inclusion functor $\inc$ restricts to a functor $\inc\colon \mathcal{F}({_{A/AeA}\Delta})\lxr \F(_{A}\Delta)$. Also, a similar argument as above shows that the inclusion functor sends $\mathcal{F}(_{A/AeA}\nabla)$ to $\mathcal{F}(_{A}\nabla)$.

(iv) This follows immediately by (ii) and (iii). 
\end{proof}
\end{cor}

\end{document}